\documentstyle[11pt,epsf]{article} 

\newtheorem{sect}{}[section]

\begin{document}

\title{Shifting Homomorphisms in Quandle Cohomology and
Skeins of Cocycle Knot Invariants}

\author{
J. Scott Carter \\
University of South Alabama \\
Mobile, AL 36688 \\ carter@mathstat.usouthal.edu \and
Daniel Jelsovsky \\
University of South Florida \\
Tampa, FL 33620 \\ jelsovsk@math.usf.edu \and
Seiichi Kamada \\
Osaka City University \\
Osaka 558-8585, JAPAN\\ kamada@sci.osaka-cu.ac.jp \\
skamada@mathstat.usouthal.edu
\and
Masahico Saito \\
University of South Florida \\
Tampa, FL 33620 \\ saito@math.usf.edu
}
\maketitle

\vspace{5mm}

\begin{abstract}
Homomorphisms on quandle cohomology groups  that raise 
the dimensions by one are studied in relation to the cocycle 
state-sum invariants of knots and knotted surfaces.
Skein relations are also studied.
\end{abstract}

\vspace{5mm}

\section{Introduction}

A quandle is a set with a 
self-distributive binary operation (defined below)
whose definition was motivated from knot theory. 
A (co)homology theory was defined in \cite{CJKLS} for quandles,
which is a modification of rack (co)homology defined in \cite{FRS2}. 
State-sum invariants using quandle cocycles as weights are 
defined \cite{CJKLS} and computed for important families
of classical knots and knotted surfaces \cite{CJKS1}.
Quandle homomorphisms and virtual knots are applied to this 
homology theory \cite{CJKS2}. The invariants were applied to study 
knots, 
for example, in detecting non-invertible
knotted surfaces \cite{CJKLS}. 

In this paper, homomorphisms on quandle cohomology groups  that raise 
the dimensions by one are studied in relation to the cocycle 
state-sum invariants of knots and knotted surfaces.
Non-triviality of such homomorphisms is proved in Section 3. 
Skein relations are studied in Section 4 for some quandles.
Preliminary material is contained in Section 2.

\newpage 

\section{Definitions of Quandle (Co)Homology and Cocycle Invariants} 
\label{def}

\begin{sect}{\bf Definition.\/}
{\rm 
A {\it quandle\/}, $X$, is a set with a binary operation
$\ast$ such that 

\noindent
(I. {\sc idempotency}) for any $a \in X$, $a \ast a=a$,

\noindent
(II. {\sc right-invertibility}) for any $a, b \in X$,
 there is a unique $c \in X$ such that $ a =
c\ast b$,
and 

\noindent 
(III. {\sc self-distributivity}) for any $a, b, c \in X$, we have
 $(a \ast b) \ast c =
(a \ast c) \ast (b \ast c)$.

A {\it rack\/} is a set with a binary operation that satisfies (II) and
(III).
Racks and quandles have been studied in, for example, 
\cite{Brieskorn},\cite{FR},\cite{Joyce},\cite{K&P}, and \cite{Matveev}.

A map $f: X \rightarrow Y$ between two quandles (resp. racks) $X,Y$
is called a quandle (resp. rack) homomorphism if 
$f(a*b)=f(a)*f(b)$ for any $a, b \in X$. 
A (quandle or rack) homomorphism is a (quandle or rack) isomorphism
if it is bijective. An isomorphism between  the same quandle (or rack) 
is an automorphism.

}\end{sect}

\begin{sect}{\bf Examples.\/}\label{quanxam} 
{\rm 
Any set $X$ with the operation $x*y=x$ for any $x,y \in X$ is 
a quandle called the {\it trivial} quandle. 
The trivial quandle of $n$ elements is denoted by $T_n$.

Any group $G$ 
is a quandle by conjugation as operation:
$a*b=b^{-1}ab$ for $a,b \in G$.  
Any subset of $G$  
 that is closed under conjugation is also a quandle.

Let $n$ be a positive integer. 
For elements  $i, j \in \{ 0, 1, \ldots , n-1 \}$, define
$i\ast j= 2j-i$ where the sum on the right is reduced mod $n$. 
Then $\ast$ defines a quandle 
structure  called the {\it dihedral quandle},
 $R_n$.
This set can be identified with  the 
set of reflections of a regular $n$-gon
 with conjugation
as the quandle operation. 
We also represent the elements of $R_3$ by  $\alpha, \beta,$ and $\gamma$, 
 where the quandle 
multiplication is given by 
$x*y = z$ where $z\ne x, y$ when $x \ne y$ and $x*x=x$, for 
$x,y,z \in \{ \alpha, \beta, \gamma \}$.

Any $\Lambda={\bf Z}[T, T^{-1}]$-module $M$ is a quandle with 
$a*b=Ta+(1-T)b$, $a,b \in M$, called an {\it  Alexander  quandle}. 
Furthermore for a positive integer 
$n$, a {\it mod-$n$ Alexander  quandle}
${\bf Z}_n[T, T^{-1}]/(h(T))$
is a quandle 
for 
a Laurent polynomial $h(T)$.
The mod-$n$ Alexander quandle is finite 
if the coefficients of the  
highest and lowest degree terms 
of $h$  
 are $\pm 1$.

See \cite{Brieskorn},
 \cite{FR}, \cite{Joyce}, or \cite{Matveev} 
 for further examples.

}\end{sect}

\begin{sect}{\bf Remark.\/}
{\rm 
Let $X$ denote a quandle.
{}From Axiom~II, each element $b \in X$ defines a bijection
$S(b) : X \to X$ with $aS(b) = a \ast b$. The bijection is an automorphism
by Axiom~III.
For a word $w = b_1^{\epsilon_1} \dots b_n^{\epsilon_n}$ 
where
$b_1, \dots, b_n
\in X;
\epsilon_1, \dots, \epsilon_n \in \{\pm 1\}$,
we define
$a \ast w = aS(w)$ by
$aS(b_1)^{\epsilon_1}\dots S(b_n)^{\epsilon_n}$.
An automorphism of $X$ is called an {\it inner-automorphism\/} of $X$ if it
is $S(w)$ for a
word $w$. (The notation $S(b)$ follows Joyce's paper \cite{Joyce} and $a
\ast w$ ($= a^w$)
follows Fenn-Rourke \cite{FR}.)
}\end{sect}

\vspace{5mm} 

 Let $C_n^{\rm R}(X)$ be the free 
abelian group generated by
$n$-tuples $(x_1, \dots, x_n)$ of elements of a quandle $X$. Define a
homomorphism
$\partial_{n}: C_{n}^{\rm R}(X) \to C_{n-1}^{\rm R}(X)$ by \begin{eqnarray}
\lefteqn{
\partial_{n}(x_1, x_2, \dots, x_n) } \nonumber \\ && =
\sum_{i=2}^{n} (-1)^{i}\left[ (x_1, x_2, \dots, x_{i-1}, x_{i+1},\dots, x_n) \right.
\nonumber \\
&&
- \left. (x_1 \ast x_i, x_2 \ast x_i, \dots, x_{i-1}\ast x_i, x_{i+1}, \dots, x_n) \right]
\end{eqnarray}
for $n \geq 2$ 
and $\partial_n=0$ for 
$n \leq 1$. 
 Then
$C_\ast^{\rm R}(X)
= \{C_n^{\rm R}(X), \partial_n \}$ is a chain complex.

Let $C_n^{\rm D}(X)$ be the subset of $C_n^{\rm R}(X)$ generated
by $n$-tuples $(x_1, \dots, x_n)$
with $x_{i}=x_{i+1}$ for some $i \in \{1, \dots,n-1\}$ if $n \geq 2$;
otherwise let $C_n^{\rm D}(X)=0$. If $X$ is a quandle, then
$\partial_n(C_n^{\rm D}(X)) \subset C_{n-1}^{\rm D}(X)$ and
$C_\ast^{\rm D}(X) = \{ C_n^{\rm D}(X), \partial_n \}$ is a sub-complex of
$C_\ast^{\rm
R}(X)$. Put $C_n^{\rm Q}(X) = C_n^{\rm R}(X)/ C_n^{\rm D}(X)$ and 
$C_\ast^{\rm Q}(X) = \{ C_n^{\rm Q}(X), \partial'_n \}$,
where $\partial'_n$ is the induced homomorphism.
Henceforth, all boundary maps will be denoted by $\partial_n$.

For an abelian group $G$, define the chain and cochain complexes
\begin{eqnarray}
C_\ast^{\rm W}(X;G) = C_\ast^{\rm W}(X) \otimes G, \quad && \partial =
\partial \otimes {\rm id}; \\ C^\ast_{\rm W}(X;G) = {\rm Hom}(C_\ast^{\rm
W}(X), G), \quad
&& \delta= {\rm Hom}(\partial, {\rm id})
\end{eqnarray}
in the usual way, where ${\rm W}$ 
 $={\rm D}$, ${\rm R}$, ${\rm Q}$.

\begin{sect}{\bf Definition \cite{CJKLS}.\/} {\rm
The $n$\/th {\it quandle homology group\/}  and the $n$\/th
{\it quandle cohomology group\/ } \cite{CJKLS} of a quandle $X$ with coefficient group $G$ are
\begin{eqnarray}
H_n^{\rm Q}(X; G) 
 = H_{n}(C_\ast^{\rm Q}(X;G)), \quad
H^n_{\rm Q}(X; G) 
 = H^{n}(C^\ast_{\rm Q}(X;G)). \end{eqnarray}

\begin{sloppypar}
The cycle and boundary groups 
(resp. cocycle and coboundary groups)
are denoted by $Z_n^{\rm Q}(X;G)$ and $B_n^{\rm Q}(X;G)$
(resp.  $Z^n_{\rm Q}(X;G)$ and $B^n_{\rm Q}(X;G)$), 
 so that
$$H_n^{\rm Q}(X;G) = Z_n^{\rm Q}(X;G)/ B_n^{\rm Q}(X;G),
\; H^n_{\rm Q}(X;G) = Z^n_{\rm Q}(X;G)/ B^n_{\rm Q}(X;G).$$
We will omit the coefficient group $G$ if $G = {\bf Z}$ as usual.
\end{sloppypar}

}
\end{sect}

\begin{figure}
\begin{center}
\mbox{
\epsfxsize=3in
\epsfbox{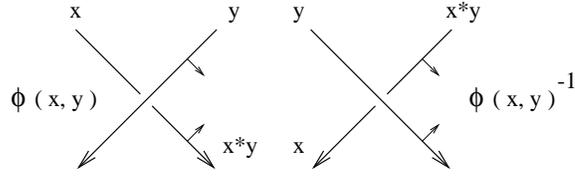}
}
\end{center}
\caption{Crossings and weights }
\label{twocrossings}
\end{figure}

\begin{figure}
\begin{center}
\mbox{
\epsfxsize=3in
\epsfbox{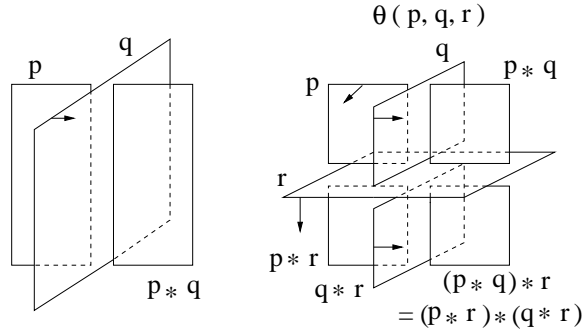}
}
\end{center}
\caption{Colors at double curves and weights at triple points }
\label{triplepoint}
\end{figure}

Assume that a finite quandle $X$ is given. 
Pick a 
quandle 
$2$-cocycle 
$\phi \in  Z_{\rm Q}^2(X; G),$ 
and write the coefficient 
group, $G$, multiplicatively. 
Consider a  crossing 
in the diagram.
For each coloring of the diagram, evaluate 
the $2$-cocycle 
 on the quandle colors that 
appear near the crossing as described as follows: 
The first argument is the color on the under-arc away from which 
the normal to the over-arc points. 
The second argument is the color on the over-arc.
Let $\tau$ denote 
a  crossing, 
let $\epsilon (\tau)$ denote its sign, 
and let ${\cal  C}$ denote a coloring.

\begin{sect}{\bf Definition \cite{CJKLS}.\/} {\rm
When the colors of 
the arcs 
 are as describe above,
the 
{\it (Boltzmann) weight
of a crossing} is  
$B(\tau, {\cal C}) = \phi(x,y)^{\epsilon (\tau)}$. 

The  {\it state-sum}, 
is the expression 
$$
\Phi_{\phi}(K) = \sum_{{\cal C}}  \prod_{\tau}  B( \tau, {\cal C}).
$$
The product is taken over all crossings of the given diagram,
and the sum is taken over all possible colorings.
The values of the 
state-sums 
are  taken to be in  the group ring ${\bf Z}[G]$.
} \end{sect}

 The coloring 
situation 
and the weights are 
depicted in Fig.~\ref{twocrossings}. 
The state-sum invariant is similarly defined for knotted surfaces in $4$-space
using coloring conventions along double curves on projections
(as in Fig.~\ref{triplepoint} left) and weights assigned to triple points
(as in Fig.~\ref{triplepoint} right).
In this case, signs $\epsilon=\pm 1$ are defined for triple points on projections, and the Boltzmann weight is defined by 
$\theta(p,q,r)^{\epsilon}$ using 
$\theta \in Z_{\rm Q}^3(X; G)$. 
 See \cite{CJKLS} for details.

It was proved in \cite{CJKLS} that for classical knots and knotted surfaces 
in ${\bf R}^4$, the state-sums are knot  invariants, by 
showing the invariance under Reidemeister moves and their $4$-dimensional
analogues (Roseman moves).

\begin{figure}
\begin{center}
\mbox{
\epsfxsize=6in
\epsfbox{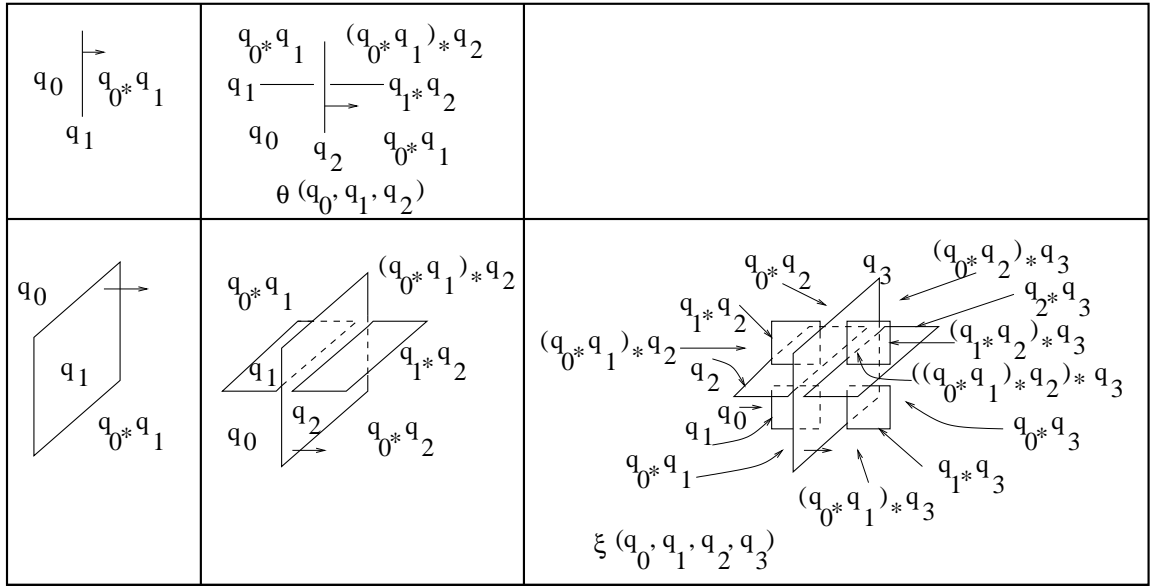}
}
\end{center}
\caption{Shadow color conventions and weights at crossings and triple points }
\label{shadowcolors}
\end{figure}

Furthermore, shadow colorings are defined using complementary regions, 
in addition to arcs and sheets, and 
 are used to define
state-sum invariants. 
Shadow colors are defined in \cite{FRS2} and used in \cite{RS}. 
The conventions and cocycle weights are depicted in Fig.~\ref{shadowcolors}.
A shadow coloring is required to satisfy the condition depicted
 in the top left entry of Fig.~\ref{shadowcolors}.
For $\theta \in Z^3_{\rm Q}(X;G)$, the value $\theta (q_0, q_1, q_2)$
corresponds to a crossing as depicted in top middle entry 
of  Fig.~\ref{shadowcolors}.
The Boltzmann weight is defined by 
$B(\tau, {\cal C})=\theta (q_0, q_1, q_2)^{\epsilon(\tau)} $
for a shadow coloring ${\cal C}$, where $\epsilon(\tau)$ is the sign of
the crossing $\tau$. Then the state-sum is defined by
$\Phi_{\theta} (K)=\sum_{{\cal C}}  \prod_{\tau}  B( \tau, {\cal C})$. 

State-sum invariants using shadow colorings are similarly defined 
using $\xi \in Z^4_{\rm Q}(X;G)$, using the coloring convention
depicted in Fig.~\ref{shadowcolors} bottom left and middle,
and using the correspondence to $\xi (q_0, q_1, q_2, q_3)$ 
depicted in Fig.~\ref{shadowcolors} bottom right.

They are also invariants of classical knots and knotted surfaces.
See \cite{CJKS2,SSS2} for more details.

\section{Shifting Homomorphisms}

\begin{sect}{\bf Definition.\/} 
{\rm 
Let $X$ be a quandle and $G$ be a coefficient abelian group.
Let $\rho=\rho^{(n)} : X^{n} \rightarrow X^{n-1}$  be 
the 
map 
defined by 
$$\rho (x_1, x_2, \cdots, x_n)  = \left\{ \begin{array}{lr}
 0 & {\mbox{\rm if }} \  x_1 = x_2 , \\
(x_2, \cdots, x_n)  & {\mbox{\rm if }} \  x_1 \neq x_2. \end{array}\right.$$
Extending linearly we obtain 
$\rho_{\sharp}=\rho_{\sharp}^{(n)}: C^{\rm R}_{n}(X; G) \rightarrow C^{\rm R}_{n-1}(X; G)$.
It is checked by computation that 
$\rho_{\sharp} (C^{\rm D}_{n}(X; G) ) \subset C^{\rm D}_{n-1}(X; G)$, so that
$\rho_{\sharp}$ induces 
$\rho_{\sharp}=\rho_{\sharp}^{(n)}: C^{\rm Q}_{n}(X; G) \rightarrow C^{\rm Q}_{n-1}(X; G)$ (the same notation $\rho_{\sharp}$ is used).
Similarly define $\rho^{\sharp}=\rho_{(n)}^{\sharp}: C^{n}_{\rm Q}(X; G) \rightarrow C^{n+1}_{\rm Q}(X; G)$ by  
$(\rho^{\sharp} f)(x)=f(\rho_{\sharp}(x))$ for all $x \in C_{n+1}^{\rm Q}(X; G)$.
} \end{sect}

By computation we have 

\begin{sect} {\bf Lemma.\/}
$\rho_{\sharp} \partial = - \partial \rho_{\sharp} $,
and $\rho^{\sharp} \delta = - \delta \rho^{\sharp} $.
\end{sect}

Hence $\rho$ induces homomorphisms on 
homology and cohomology groups, 
denoted by $\rho_*$ and $\rho^*$.
We call the
homomorphisms 
$\rho_\#$,  $\rho^\#$,  $\rho_*$,  $\rho^*$
the {\it shifting homomorphisms.}

\begin{sect} {\bf Proposition.\/} 
 Let $\phi \in  C^{n}_{\rm Q}(X; G)$ 
be an $n$-cochain for some
quandle $X$ and an abelian group $G$,  and $\rho$ be as above.  
 Then $\rho^{\sharp} \phi$ is an
$(n+1)$-cocycle (i.e. 
$\rho^{\sharp} \phi \in Z^{n+1}_{\rm Q} (X; G)$)
if and only if $\phi$ is an $n$-cocycle ($\phi \in Z^n_{\rm Q}(X; G)$).
\end{sect}
{\it Proof.\/}
We have  that $\delta \rho^\sharp = -\rho^\sharp \delta$.
 So if $\phi \in Z^{n}_{\rm Q}(X; G)$, 
then 
$$
\delta \rho^\sharp \phi (x_0, \ldots , x_{n}) 
= - \rho^\sharp \delta \phi (x_0, \ldots , x_{n}) = 0. 
$$
If  $\rho^\sharp \phi \in Z^{n+1}_{\rm Q}(X;G) $, then 
$\delta \rho^\sharp \phi (x_0, \ldots, x_{n}) =0$.
On the other hand, 
$\delta \rho^\sharp \phi (x_0, \ldots, x_{n}) =\delta \phi (x_1, \ldots , x_{n} )$ if $x_0 \neq x_1$, so the result follows. $\Box$

\begin{sect}{\bf Remark.\/}
{\rm For some positive integer $n$,
consider $\rho_{\sharp} : C_{n}^{\rm Q}(X) \rightarrow C_{n-1}^{\rm Q}(X)$.
 Clearly if
$\eta \in Z_n^{\rm Q}(X)$  
is a cycle, then
$\partial_{n-1}\rho_\sharp \eta=-\rho_\sharp \partial_n \eta =0$.
However, if 
$\rho_{\sharp}(\zeta)$ 
is a cycle, 
the  chain $\zeta$ need not be a cycle.
For example, let  $\zeta=(2,1,3) \in C_3(R_4)$, 
we see
that $\rho          (2,1,3)=(1,3)$ and $\partial_2(1,3)=1-1*3=1-1=0$.
However, $\partial_3(2,1,3) = (2,3)-(0,3)-(2,1)+(0,1) \neq 0$.
Thus, $(2,1,3)$ is not a 3-cycle.
} \end{sect}

\begin{figure}
\begin{center}
\mbox{
\epsfxsize=4.5in
\epsfbox{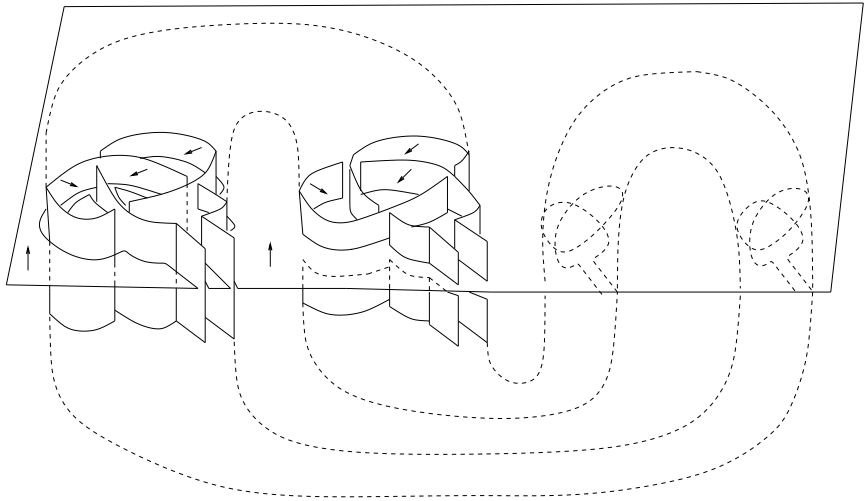}
}
\end{center}
\caption{A diagram of $2$-twist spun trefoil }
\label{shin}
\end{figure}

\begin{sect} 
{\bf Proposition.\/}
$ H^4_{\rm Q}(R_3; {\bf Z}_3) \neq 0 $.
 \end{sect}
{\it Proof.\/}
It was proved in \cite{CJKLS} that if a coboundary is used to 
define the  state-sum invariant, then the invariant is trivial 
(i.e., every state-sum term is $1$).
Analogous arguments  apply to the state-sum 
invariants defined using shadow colorings to prove that 
 if $\xi$ is a coboundary, then 
the state-sum term is trivial (c.f. \cite{SSS2}). 
Hence we prove that there exists a cocycle
 $\xi \in Z^4_{\rm Q}(R_3; {\bf Z}_3)$ and a shadow coloring of a knotted surface diagram with a non-trivial state-sum term.

Specifically, let $K$ be the $2$-twist spun trefoil, and 
let $\phi \in   Z^3_{\rm Q}(R_3; {\bf Z}_3)$ be a $3$-cocycle 
defined (in \cite{CJKLS}) by 
$$ \phi= \chi_{0,2,0} + \chi_{1,0,1} 
+ \chi_{1,0,2} + \chi_{2,0,2} + \chi_{2,0,2} + \chi_{2,1,2}
- \chi_{0,1,0} - \chi_{0,2,1}, $$
where  $\chi$ denotes the characterictic function
$$\chi_{a,b, c} (x,y, z) = \left\{ \begin{array}{lr} 1 & {\mbox{\rm if }} \ 
(x,y, z)=(a,b, c), \\
0 & {\mbox{\rm if }} \ 
(x,y,z)\not=(a,b,c). \end{array}\right.$$
We show that there is a shadow coloring of a diagram of
$K$ which gives a non-trivial
state-sum term with the cocycle 
 $\xi= \rho^{\sharp} \phi  \in Z^4_{\rm Q}(R_3; {\bf Z}_3)$. 

 For this purpose we use a diagram of $K$ given in \cite{Shin},
as depicted in Fig.~\ref{shin}. 
In the figure, a diagram of trefoil with a small portion removed 
goes around the horizontal plane twice. 
The second time around (the right side of the figure)
is depicted only schematically.
The left-most trefoil goes over the plane with respect to the height 
direction, and the second goes 
under the plane, and  this pattern is repeated one more time on the 
right side of the figure.
The portion removed are connected to
the plane by branch points, which are not depicted, as
these portions of the diagram do not contain triple points and hence
do not make any contribution to the state-sum invariant.
The plane  in fact represents a sphere. The continuous trace of 
the moving trefoil, together with the sphere attached to it via branch
points, gives a diagram of $2$-twist spun trefoil.
The orientation normal vectors are also depicted in the figure.
In \cite{Shin} it was shown that with this diagram the state-sum 
is computed as 
\begin{eqnarray*} 
\lefteqn{ \Phi_{\phi} (K) = 
 3+  2 \ [ \ \phi(0,2,1)^{-1} \ \phi(2,1,0) \ \phi(0,1,2)^{-1} \ \phi(1,2,0) } \\
& & +  \phi(1,0,2)^{-1} \ \phi(0,2,1) \ \phi(1,2,0)^{-1} \ \phi(2,0,1) \\
& & +  \phi(2,1,0)^{-1} \ \phi(1,0,2) \ \phi(2,0,1)^{-1} \ \phi(0,1,2) \  ]. 
\end{eqnarray*}
Each of these terms containing $\phi$
contributes $t$ (a generator of ${\bf Z}_3$ 
in multiplicative notation), giving the state-sum invariant $3+6t$. 

Consider the region $R$ of the diagram which lies below the horizontal plane
and is outside of the trace of the trefoil diagram. 
{}From the normal vectors, it is seen that the color of $R$ appears in 
the first entry of a $4$-cocycle $\xi$, when the state-sum invariant is 
computed using  $\xi$ and shadow colorings. 
Note also that given a coloring of the above diagram, 
any choice of color for $R$ extends to a shadow coloring of the diagram
(c.f. \cite{SSS2}). 
Hence we choose a shadow coloring with $2 \in R_3$ as 
the color assigned to $R$. With this choice, there are three terms
in the above expression of $\phi$ that give the repetitive first and second 
entries 
 for $\xi$: $\xi(2,2,1,0)$, $\xi(2, 2,0,1)$, and $\xi(2,2,1,0)$. 
However, the corresponding triples  evaluates trivially by $\phi$.
Therefore, this shadow coloring contributes $t$ to the 
invariant. The result follows.
$\Box$

\begin{sect} 
{\bf Remark.\/} {\rm 
It is an interesting problem to determine when
 $\rho^{\sharp}_{(n)}$ is injective.
The above proof illustrates the use of colored knot diagrams
and shifting homomorphisms
for solving this problem.
 In particular, 
colored knot diagrams are defined for higher dimensions \cite{SSS2}, 
and the above method can be stated as follows as a conjecture.

 {\it Conjecture\/: Let $X$ be a quandle, $G$ an abelian 
group, $\phi \in Z^n_{\rm Q} (X; G)$, and $K$ 
be an $(n-1)$-knot diagram 
in  ${\bf R}^n$.
If there is a color ${\cal C}$ of $K$ 
with a non-trivial state-sum term for $\Phi_{\phi}(K)$, 
then  $H^{n+1}_{\rm Q}(X;G) \neq 0.$ } 
} \end{sect}


\begin{sect} 
{\bf Proposition.\/}
\label{surj} Let $k$ and $m$ be positive
integers.
Let $X$ be the Alexander quandle 
$X={\bf Z}_{mk}[T,T^{-1}]/(T-1 \pm m)$.
The orbit quandle ${\rm Orb}(X) $ is isomorphic to the trivial quandle $T_m$.
\end{sect}
{\it Proof.\/}  
We represent the elements of $X$ as $\{0,1, \ldots , mk-1\}$ with
quandle multiplication given by $a*b= (1\mp m)a \pm m b {\pmod{ mk}}$.
Consider $f : X \rightarrow T_m$ given by $f(a)=a$ (mod $m$). 
 Then 
$$f(a*b)    
= (1 \mp m)a \pm mb 
=a 
=f(a) * f(b)
\pmod{m}.  
$$
Thus $f:X\rightarrow T_m$ is a  surjective 
quandle homomorphism 
which induces a surjective  quandle homorphism $\tilde{f}:{\rm Orb}(X) \rightarrow T_m.$ Suppose that $f(a)= f(b)$, so $b=a+ms$ for some integer $s$.
 Then $a*(a\pm s)=b$ 
(for $X={\bf Z}_{mk}[T,T^{-1}]/(T-1 \pm m)$ respectively), 
and $a$ and $b$ are in the same orbit. Thus 
$\tilde f$ is injective. $\Box$

\begin{sect} 
{\bf Corollary.\/}
${\rm Orb} ({\bf
Z}_{8}[T,T^{-1}]/(T-3)) =T_2, \quad  {\rm Orb} ({\bf
Z}_{8}[T,T^{-1}]/(T-5)) =T_4$.
\end{sect}

Observe that either result can be obtained by direct computation since the 
quandles in question have $8$ elements.
In \cite{CJKS2}  it was proved that 
${\rm rank}(H^n_{\rm Q}(X;G)) \geq {\rm rank}(H^n_{\rm Q}({\rm Orb}(X);G))$
by considering the pull-backs of elements of $C^n_{\rm Q}({\rm Orb}(X);G)$
in $C^n_{\rm Q}(X;G)$.   Using the shifting homomorphism, we
have the following.
\begin{sect} {\bf  Theorem.\/}
For $X= {\rm Z}_8[T,T^{-1}]/(T-5)$,
there is 
a non-coboundary 
$n$-cocycle $\in Z^n_{\rm Q}( X; {\bf Z}_2)$
which is not the pull-back of an $n$-cycle
$\in Z^n_{\rm Q}( T_4; {\bf Z}_2)$ for $n=2,3$.
As a consequence, 
${\rm rank}(H^n_{\rm Q}(X;G)) > {\rm rank}(H^n_{\rm Q}(T_4;G))$.
\end{sect}  
{\it Proof.\/} For $H^2$, 
we used {\sc Maple} to calculate that
the value $\Phi_\theta(T(2,4))=48 +16 t$
 where  the 
cocycle
$\theta =  {\chi _{0,
\,1}} + {\chi _{0, \,5}} + {\chi _{1, \,5}} + {\chi _{ 2, \,1}} +
{\chi _{2, \,5}} + {\chi _{3, \,5}} + {\chi _{5, \,1}}
 + {\chi _{7, \,1}}$
takes values in ${\bf Z}_2$,
and $T(2,4)$ denotes the $(2,4)$-torus link.

  Now, consider a cocycle $\phi \in Z^2_Q(X,{\bf Z}_2)$
which is the pull back of a $2$-cocycle in
 $T_4$.  Then $$\phi = \sum_{i,j \in \{0,1,2,3\}, 
i \neq j} a_{i,j}
 (\chi_{i,j} + \chi_{i+4,j} + \chi_{i,j+4} + \chi_{i+4,j+4}),$$
 where the $a_{i,j}\;$s 
 are constants.  For any term $a_{i,j}
 (\chi_{i,j} + \chi_{i+4,j} + \chi_{i,j+4} + \chi_{i+4,j+4})$,
 consider its state-sum contribution for $T(2,4)$ with the coefficient group
 ${\bf Z}_2$.  We may
 consider $T(2,4)$ as the closure of the braid $\sigma_1^4$.  Then
 the only time any crossing will have non-trivial weight with
 respect to this particular $a_{i,j}
 (\chi_{i,j} + \chi_{i+4,j} + \chi_{i,j+4} + \chi_{i+4,j+4})$ is
 when the initial color vector of the braid is of the form $(x,y)$
 where $x \equiv i$ (mod 4) and $y \equiv j$ (mod 4).  However, if
 this is the case, we have two crossings of this form.  So the
 state-sum contribution of this color is still trivial (since
 $t^2=1$ in ${\bf Z}_2$).  Thus, this part of the cocycle
 contributes $1$ 
 to the state-sum.  Since the entire cocycle
 is made up of these parts, the state-sum of $T(2,4)$ with  ${\bf
 Z}_2$ is trivial, and, in fact, is $64$.  Since $64 \neq 48+16t$,
 $\theta$ and $\phi$ are not cohomologous.  So, $\theta$ is not a
 pulled back cocycle, which means it is not accounted for by the
 rank of the trivial quandle.


\begin{figure}
\begin{center}
\mbox{
\epsfxsize=.8in
\epsfbox{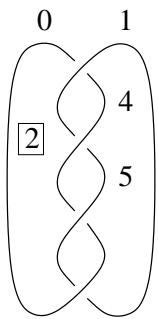}
}
\end{center}
\caption{A shadow coloring of $T(2,4)$ }
\label{shadowT24}
\end{figure}

For $H^3$, 
 the state-sum invariant of $K$, using $\rho^{\sharp} \theta$
(where $\theta$ is as above) and using  shadow colorings of $T(2,4)$,
is non-trivial.
This can be seen using the shadow coloring depicted in Fig.~\ref{shadowT24}.

 If $\phi$
is a pullback of a $3$-cocycle in 
$Z^3_Q(T_4;{\bf Z}_2)$, then it
is of the form 
\begin{eqnarray*} \sum_{i,j,k \in {0,1,2,3},\/ i \neq j,\/ j \neq
k} a_{i,j,k}
 (\chi_{i,j,k} 
 + \chi_{i+4,j,k} + \chi_{i,j+4,k} + \chi_{i+4,j+4,k}
& & \\
 + \chi_{i,j,k+4} + \chi_{i+4,j,k+4} + \chi_{i,j+4,k+4} +
\chi_{i+4,j+4,k+4}). \end{eqnarray*} 
For any term  $$a_{i,j,k}
 (\chi_{i,j,k} + \chi_{i+4,j,k} + \chi_{i,j+4,k} + \chi_{i+4,j+4,k}+
 \chi_{i,j,k+4} + \chi_{i+4,j,k+4} + \chi_{i,j+4,k+4} +
 \chi_{i+4,j+4,k+4})$$
consider its state-sum contribution for $K$ over ${\bf Z}_2$.
First note that for a given color the first term of every weight
is the same.  This is since
the color on the region at infinity
 remains the same
for all crossings of $T(2,4)$. Then, as with the $2$-dimensional
case, every choice of colors will either be trivial at every
crossing, or cancel each other out.  Thus, the state-sum invariant
of $K$ for $\phi$ is $512$, and hence $\rho^{\sharp} \theta$ is not a
pullback from the trivial quandle.
$\Box$

\section{Skein Relations}

Being state-sum invariants, the cocycle invariants
have skein relations coming from the minimal polynomials
of the corresponding $R$-matrices.
However, the $R$-matrices are in general too large to
compute by computers. 
In this section we give a method of finding skein relations 
using Burau matrices, which are much smaller than $R$-matrices. 
By skein relations, we mean linear formulas for 
the state-sum invariant, and we do not study the sufficiency of 
these as recursive formulas (i.e., whether or not they form 
a complete set of formulas to compute the invariant 
for any knots and links recursively).

The coefficient groups  of cocycle groups used in this section  
are cyclic groups $G={\bf Z}_p$ for some integer $p$, that 
are denoted multiplicatively, 
 ${\bf Z}_p= \langle t | t^p=1 \rangle= \{ t^n | n=0, 1, \cdots, p-1 \}$.
In this case, the state-sum takes values in 
${\bf Z}[G]={\bf Z}[t]/(t^p-1)$.

An oriented $n$-tangle for a positive integer $n$ is 
a tangle (diagram) with $n$ strings going in at the top, and with $n$ 
strings going out from the bottom of the diagram. 
Suppose that  the colors on the top and 
bottom strings of a tangle are specified  by vectors 
 $[c_1, \cdots, c_n]$ and $[c'_1, \cdots, c'_n]$,
respectively, such that  the entries are 
elements of  $X$.
The state-sum for such a tangle is defined similarly,
and denoted by 
$$\Phi(T) \left[ \begin{array}{lll} c_1, &  \cdots,& c_n
\\ c'_1, &  \cdots, & c'_n \end{array} \right]
= \sum_{\cal C} \prod_{\tau} B(\tau, {\cal C}), $$
where $ {\cal C}$ ranges over all colorings that restrict to 
the given colors, $[c_1, \cdots, c_n]$ and $[c'_1, \cdots, c'_n]$, 
on the boundary (top and bottom) segments, and $\tau$ ranges over all 
crossings of the tangle $T$.

\begin{sect} {\bf Lemma.\/} \label{skeinlem}
Let $T_i$, $i=1, \cdots, m$, be $n$-tangles for 
positive 
integers $m,n$. 
Let $K_i$,  $i=1, \cdots, m$, be classical 
knot or link diagrams 
such that they are all identical outside of a small ball inside which 
they have $n$-tangles $T_i$ respectively. 
Suppose that 
\begin{itemize}
\item 
the set of color vectors on the top and bottom strings of $T_i$
are identical for all $i$, in the sense that 
if $[c_1, \cdots, c_n]$ and $[c'_1, \cdots, c'_n]$
are quandle vectors that color the top and bottom (respectively) strings 
of $T_j$ for some $j$, $j=1, \cdots, n$, then they color the strings
of $T_i$ for all $i=1, \cdots, n$,
uniquely. 
\item 
There exists a set of Laurent polynomials 
$f_i(t)$, $i=1, \cdots, n$ such that 
for any top and bottom color vectors 
 $[c_1, \cdots, c_n]$ and $[c'_1, \cdots, c'_n]$
that color $T_i$s, 
the state-sum term
 $\Phi(T_i)=\Phi(T_i) \left[ \begin{array}{lll} c_1, &  \cdots,& c_n
\\ c'_1, &  \cdots, & c'_n \end{array} \right]$  

 for $T_i$
satisfy the equality 
$$ f_1(t) \Phi(T_1) + \cdots + f_n(t) \Phi(T_n) = 0.$$
\end{itemize}

Then the cocycle invariant satisfies the skein relation
$$ f_1(t) \Phi(K_1) + \cdots + f_n(t) \Phi(K_n) = 0.$$
\end{sect}
{\it Proof.\/} It follows from 
the state-sum definition 
(or by using the $R$-matrix description).
$\Box$

\begin{sect} {\bf Example.\/} {\rm
We illustrate 
the
approach of finding skein relations 
using the above lemma and Burau matrices 
over 
$R_4$ (which is the Alexander quandle ${\bf Z}_2[T,T^{-1}]/(T^2+1)$)
with the cocycle
$\phi=\chi_{(0,1)} \chi_{(0,3)}$  
(see \cite{CJKLS}). 
Let $T_+$, $T_0$, and $T_-$ be the $2$-tangles represented by 
braid words $\sigma_1^4$, $1$, and $\sigma_1^{-4}$
respectively. Then in \cite{CJKLS} it is observed that 
they satisfy the first condition of Lemma~\ref{skeinlem}.
For the color vectors $[i,i]$, $[2i, 2j]$, or $[2i+1, 2j+1]$
 on top strings, for any $i,j$,
the state-sum contribution is trivial ($1\in G$) 
for all $T_k$. 
Thus the skein expression 
$$ f_+(t) \Phi(T_+) +  f_0(t) \Phi(T_0) + f_-(t) \Phi(T_-) = 0 $$ 
gives 
$$ f_+(t)  +  f_0(t)  + f_-(t)  = 0 $$
for Laurent 
polynomials $f_k(t)$ that are to be determined. 
For other color vectors, $\Phi(T_+)=t$, $\Phi(T_0)=1$, and 
$\Phi(T_-)=t^{-1}$. Hence we have 
$$ f_+(t) \; t  +  f_0(t)  + f_-(t)  \; t^{-1}  = 0. $$
The 
choices $f_+(t)=1-t^{-1}$, $f_0=t^{-1}+t$, and $f_-(t)=t-1$
gives a solution, and we obtain a skein relation 
$$(1-t^{-1}) \Phi(K_+) - (1-t) \Phi(K_-) = (t-t^{-1}) \Phi(K_0)$$
where $K_k$ represent links with $2$-tangles $T_k$ in them.
} \end{sect}

\begin{figure}
\begin{center}
\mbox{
\epsfxsize=2.5in
\epsfbox{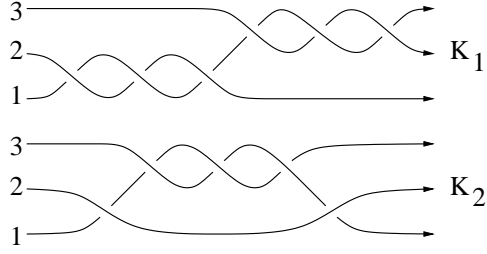}
}
\end{center}
\caption{A skein relation for $S_4$ }
\label{s4skein}
\end{figure}

\begin{sect} {\bf Example.\/} \label{s4ex} {\rm 
In this example we examine  skein relations for $S_4={\bf Z}_2[T, T^{-1}]/(T^2+T+1)$ 
and the cocycle $ \phi= \prod_{x, y \neq T, \  x \neq y } 
\chi_{(x,y)}$. 
Here, $S_4=\{ 0,1,T, 1+T\}$, and the above product ranges over all 
pairs having no $T$, and $\chi$ denotes the characterictic function
$$\chi_a(b) = \left\{ \begin{array}{lr} t & {\mbox{\rm if}} \  \ a=b
\\
                                               1   & {\mbox{\rm if}} \ \  a\ne
b \end{array} \right.$$
for pairs $a$ and $b$. 
Here the values in 
$G={\bf Z}_2$ of cocycles are denoted multiplicatively.

First, recall from \cite{CJKS1} that the trefoil and its mirror image
has the same value $4(1+3t)$ with $S_4$ with the cocycle $\phi$.
Therefore $S_4$ has the skein relation 
$$ \Phi(K_{3+})= \Phi(K_{3-}) $$ 
where $K_{3+}$ and $K_{3-}$ denote the links with the braid word
$\sigma_1^3$ and $\sigma_1^{-3}$ in $B_2$ respectively
(and the outside of these braid words are identical as usual
with the skein relations).

Second, it is seen that the braid words depicted in Fig.~\ref{s4skein}
have the identity Burau representation with $S_4$ 
(as does $\sigma_1^3$). 
Let $K_i$, $i=1,2$ be links with the braids depicted in the figure,
and let $K_0$ be the link with the identity braid word in place. 
Then we set up the skein relation 
$$ f_0(t) \Phi(K_0) +  f_1(t) \Phi(K_1) + f_2(t) \Phi(K_2) = 0 .$$
Give the numbers $1$, $2$, and $3$ on the left arcs of
 Fig.~\ref{s4skein} as depicted, from bottom to top.
Let $C_i$ be the colors on these arcs, for $i=1,2,3$.
We have the following cases. 

Case $(A)$: $C_1=C_2=C_3$ or $C_i$ all distinct and $C_1 * C_2=C_3$.
In this case both  $K_1$ and $K_2$ contribute $1$  to the state-sum 
and hence gives a relation 
$$ f_0 + f_1 + f_2=0. $$

Case $(B)$: $C_1=C_2 \neq C_3$ or $C_1 \neq C_2 = C_3$. 
In this case both  $K_1$ and $K_2$ contribute $t$  to the state-sum 
and hence gives a relation 
$$ f_0 + t (f_1 + f_2)=0. $$

Case $(C)$: $C_1=C_3 \neq C_2$ or $C_i$ all distinct and $C_1 * C_2 \neq C_3$. 
In this case $K_1$ contributes $1$ and $K_2$ contributes $t$.
Therefore we obtain
$$ f_0 + f_1 +  t f_2=0. $$

The three conditions reduce to 
$(t-1)f_i=0$ for $i=1,2$ and $f_0=-(f_1+f_2)$.
For example we obtain the relation
$$ (t+1) ( \Phi(K_1) + \Phi(K_2) - 2 \Phi(K_3) ) = 0.  $$
The factor $(t+1)$ cannot be removed from this relation because it is 
not a unit. 

Third, {\sc Maple} 
computations show that with $S_4$ the Burau matrices 
coincide for the following variations of braid words
(corresponding to the figure-eight knot), and 
the state-sum term also coincide for every color on the top strings.
Hence the invariant does not change by replacement of one by another
among the following: 
$\sigma_1 \sigma_{2}^{-1} \sigma_1 \sigma_{2}^{-1}$,
$\sigma_{1}^{-1} \sigma_{2} \sigma_{1}^{-1} \sigma_{2}$,
$\sigma_{2} \sigma_{1}^{-1} \sigma_{2} \sigma_{1}^{-1}$,
$\sigma_{2}^{-1} \sigma_{1} \sigma_{2}^{-1} \sigma_{1}$.

} \end{sect}

\begin{figure}
\begin{center}
\mbox{
\epsfxsize=2.5in
\epsfbox{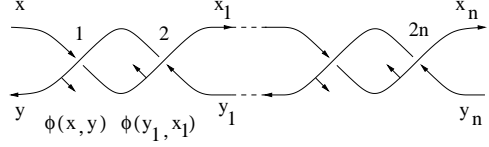}
}
\end{center}
\caption{Antiparallel strings with $2n$ crossings  }
\label{antipara}
\end{figure}

\begin{sect} {\bf Lemma.\/}
Let $P_n$ denote the  two strings with $2n$ crossings and 
with oposite orientations as depicted in Fig.~\ref{antipara}.
(If $n$ is positive, the crossings are all positive, and
for negative $n$, the 
crossings 
are understood to be negative.)

If the left two end points receive
the colors $x$ and $y$ respectively
in an Alexander quandle $X$, as depicted, then the right end points 
receive the colors $x_n$ and $y_n$ where

\begin{eqnarray*}
x_n & = & ( nT - (n-1) )x + n(1-T)y, \\
y_n & = & n(T-1)x + ((n+1)-nT)y. 
\end{eqnarray*}

In particular, if the Alexander quandle has 
coefficients in ${\bf Z}_n$,
then $x_n=x$ and $y_n=y$. 

Furthermore, the crossings contribute the following terms to the 
state-sum expression if $P_n$ is a part of a link:

$$ \prod_{i=1}^n \phi (x_i, y_i) \phi (y_{i+1}, x_{i+1}). $$

\end{sect}
{\it Proof.\/} Induction on $n$. $\Box$

\begin{sect} {\bf Proposition.\/}
If $K$ and $K'$ are related by a sequence of replacements
of $P_2$ by $P_0=I$ (the two anti-parallel strings with no crossing)
or vice versa, 
then with the cocycle $\phi \in Z^2_{\rm Q}(S_4; {\bf Z}_2)$
in Example~\ref{s4ex}, $\Phi(K)=\Phi(K')$. 
\end{sect}
 {\it Proof.\/} 
It is computed using the preceding lemma 
that every color  of $P_{2}$ contributes  $1$ to the state-sum. 
$\Box$

State-sums with shadow colors have skein relation as well.
Here we give an example with $R_3$. We label the elements 
of $R_3$ as $0,1,2$ with the 
quandle operation 
$i*j=2j-i \  \ {\pmod{3}}$.

\begin{figure}
\begin{center}
\mbox{
\epsfxsize=3in
\epsfbox{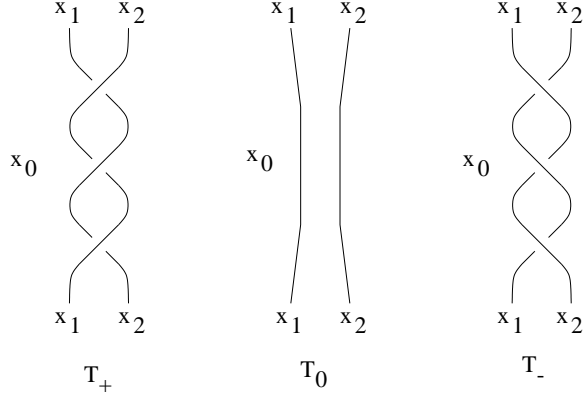}
}
\end{center}
\caption{Tangles involved in the skein relation for $R_3$ }
\label{trefskein}
\end{figure}

Let $\xi \in Z^3_{\rm Q} (R_3; {\bf Z}_3)$ denote the 
$3$-cocycle 
$$\xi = \chi_{012} \chi_{021} 
\chi_{101} \chi_{201} \chi_{202}\chi_{102}$$
where 
$$\chi_{abc} (x,y,z) = \left\{ \begin{array}{lr} t & {\mbox{\rm if }} \ 
(x,y,z)=(a,b,c), \\
1 & {\mbox{\rm if }} \ 
(x,y,z)\not=(a,b,c). \end{array}\right.$$
Let $\Phi_\xi$ denote the state-sum invariant of a link
associated with the cocycle $\xi$.  Thus 
$$\Phi_\xi  
= \sum_{\mbox{\rm shadow colorings}} \prod_{i} \xi (a_i,b_i,c_i)^{\epsilon_i}$$ 
where $(a_i,b_i,c_i)$ 
are the incoming colors at a crossing, $\epsilon_i$ is the sign of the crossing and the product ranges over all crossings.
Consider the tangle $T_+$, $T_0$, and $T_-$ that are depicted in Fig.~\ref{trefskein}.
Let $K_+, K_0, K_-$ be links with 
$T_+, T_0, T_-$ in them, respectively.

\begin{sect}{\bf Theorem.\/}
In the notation above,
$$(1-t^{-1}) \Phi_\xi(K_+) - (1-t) \Phi_\xi(K_-) = (t-t^{-1}) \Phi_\xi(K_0).$$
\end{sect}
 {\it Proof.} In the figure, let $x_0$, $x_1$, and $x_2$ denote 
colors by $R_3$ that are indicated in the figure.  Thus $x_0$ is the color on the region to the left of the tangle, $x_1$ is the color on the top left string and $x_1$ is the color on the right string.

We find
that the triple $f_{+}(t) = (1-t^{-1})$, $f_0(t)=-(t-t^{-1})$, and $f_{-}(t) = t-1$ is a solution of the equation
$$f_{+}(t) \Phi_\xi (T_+) + f_{0}(t)  \Phi_\xi (T_0)
+ f_{-}(t) \Phi_\xi (T_-)=0.$$

If $x_1 =0$ and $x_2= 1$, 
then the colors on the crossings of the  
tangles $T_+$, $T_0$ and $T_-$  are  read  in the table below.

\begin{center} 
\begin{tabular}{||c|c|c||}
\hline \hline $T_+$ & $T_0$ & $T_-$ \\ \hline \hline
$(x_0,0,1)$ & $\emptyset$ & $(x_0,2,0)^{-1}$ \\ \hline
$(x_0,1,2)$ & $\emptyset$ & $(x_0, 1,2)^{-1}$ \\ \hline
$(x_0, 2, 0)$ & $\emptyset$ & $(x_0,0,1)^{-1}$ \\ \hline \hline
\end{tabular}
\end{center} 

For any value of $x_0$, the $3$-cocycle $\xi$ evaluates to
$t$ on $T_+$, $1$ on $T_0$ and $t^{-1}$ on $T_-$. 
Thus
$$t f_+(t) + f_0(t) + t^{-1} f_-(t) =0.$$
When $x_1 \ne x_2$, then cocycle $\xi$ evaluates similarly. When
$x_1=x_2$, we obtain,
$f_+(t) + f_0(t) +f_{-}(t) =0.$  
The result follows. $\Box$

\begin{sect}{\bf Remark.\/} 
{\rm Let $\Phi'$ denote the invariant associated with
$\xi^{-1}$.  
Then 
$$(1-t)\Phi'(K_+)- (1-t^{-1}) \Phi'(K_{-}) = (t^{-1} -t) \Phi'(K_0).$$
}\end{sect}

\noindent
{\bf Acknowledgements.} JSC is being supported by NSF grant DMS-9988107.
MS is being supported by NSF grant DMS-9988101. 
SK is being supported by a Fellowship from the Japan Society for the Promotion of Science. 
We have had productive conversations with 
Shin Satoh and 
Dan Silver about this paper.


\begin{thebibliography}{99}



\bibitem{Brieskorn} Brieskorn, E., 
{\it Automorphic sets and singularities,}
Contemporary math., 78 (1988), 45--115.


\bibitem{CJKLS}
 Carter, J.S.; Jelsovsky, D.; Kamada, S.; Langford, L.; Saito, M.,
{\it Quandle cohomology and state-sum invariants
of knotted curves and surfaces,}
preprint at 
http://xxx.lanl.gov/abs/math.GT/9903135 . 

\bibitem{CJKS1}
 Carter, J.S.; Jelsovsky, D.; Kamada, S.; Saito, M.,
{\it Computations of quandle cocycle invariants of
 knotted curves and surfaces,}
preprint at 
http://xxx.lanl.gov/abs/math.GT/9906115 .

\bibitem{CJKS2}
 Carter, J.S.; Jelsovsky, D.; Kamada, S.; Saito, M.,
{\it Quandle homology groups, their betti numbers, and virtual knots,}
preprint at 
http://xxx.lanl.gov/abs/math.GT/9909161 . 
to appear in J. of Pure and Applied Algebra.

\bibitem{SSS2}
 Carter, J.S.; Kamada, S.; Saito, M.,
{\it Geometric interpretations of quandle homology,}
preprint at 
http://xxx.lanl.gov/abs/math.GT/0006115 .



\bibitem{CS:book} Carter, J.S.; Saito, M.,
{\it Knotted surfaces and their diagrams,}
the American Mathematical Society,  1998.




\bibitem{CS:cancel} Carter, J.S.; Saito, M.,
{\it Canceling branch points on the projections of surfaces in 4-space},
Proc. AMS  
116, 1, (1992) 229-237.




\bibitem{FR}   Fenn, R.; Rourke,  C.,
\textit{Racks and links in codimension two,}
Journal of Knot Theory and Its Ramifications Vol. 1 No. 4 (1992), 343-406.



 \bibitem{FRS1} 
Fenn, R.; Rourke, C.; Sanderson, B., {\it Trunks and classifying spaces,}
Appl. Categ. Structures 3 (1995), no. 4, 321--356.


\bibitem{FRS2} Fenn, R.; Rourke, C.; Sanderson, B.,
{\it James bundles and applications,} preprint found at
http://www.maths.warwick.ac.uk/${}^{\sim}$bjs/ .

\bibitem{Flower} Flower, Jean, {\it Cyclic Bordism and Rack Spaces,} Ph.D. Dissertation, Warwick (1995).

\bibitem{FoxTrip} Fox, R.H., 
{\it A quick trip through knot theory,}
in Topology of $3$-Manifolds, 
Ed. M.K. Fort Jr., Prentice-Hall (1962) 120--167.



\bibitem{Greene} Greene, M. T. {\it Some Results in Geometric Topology and Geometry,} Ph.D. Dissertation, Warwick (1997).



\bibitem{Joyce} Joyce, D.,
{\it A classifying invariant of knots, the knot quandle,}
J. Pure Appl. Alg., 23, 37--65.


\bibitem{NK} Kamada, N., {\it Alternating link diagrams on compact oriented surfaces}, (1995) preprint.


\bibitem{KK} Kamada, N. and Kamada S., 
{\it Abstract Link Diagrams and Virtual Knots}, Journal of Knot
Theory and its Ramifications 9 No.1 (2000), 93-106.




\bibitem{Lou}  Kauffman,  L. H.,
{\it Virtual knots,} Europ. J. Combinatorics (1999) 20, 663-690.

\bibitem{K&P}  Kauffman, L. H., {\it Knots and Physics},
World Scientific, Series on knots and everything, vol. 1, 1991.


\bibitem{Matveev} Matveev, S., 
{\it Distributive groupoids in knot theory,} (Russian) Mat. Sb. (N.S.)
119(161) (1982), no. 1, 78--88, 160.



\bibitem{Rose} Roseman, D.,
{\it Reidemeister-type moves for surfaces in four dimensional space, } 
in Banach Center Publications 42 (1998) Knot theory, 347--380.

\bibitem{Rosi} Rosicki, Witold, 
{\it Some simple invariants of the position of a surface in ${\bf R}^4$}, 
     Bull.of the Pol. Ac.of Sci. Math. 46(4), 1998, 335-344.
                                 

\bibitem{RS} Rourke, C., and Sanderson, B.,
{\it There are two $2$-twist-spun trefoils,}
preprint at 
 http://xxx.lanl.gov/abs/math.GT/0006062 .

\bibitem{Shin} Satoh, S., 
{\it Surface diagrams of twist-spun $2$-knots,}
preprint.


\end{thebibliography}
\end{document}